\documentclass[a4paper]{article}
\usepackage[dvips]{graphicx}
\usepackage{amssymb}
\usepackage{amsmath}
\usepackage{color}

\newcommand{\C}{\mathbb{C}}

\makeatletter
\def\dddots{\mathinner{\mkern1mu\raise\p@
    \hbox{.}\mkern2mu\raise4\p@\hbox{.}\mkern2mu
    \raise7\p@\vbox{\kern7\p@\hbox{.}}\mkern1mu}}%
\makeatother
\newtheorem{theorem}{Theorem}
\newtheorem{lemma}{Lemma}

\newtheorem{prop}{Proposition}
\newtheorem{defn}{Definition}

\newfont{\BBB}{msbm10 scaled \magstep1}
\newfont{\BBS}{msbm10}

\usepackage{mathtools,ulem} 
\usepackage[nodisplayskipstretch]{setspace}

\newcommand{\acnst}[1]{{\color{red}#1}}

\newcommand{\mathd}{\mathrm{d}}
\newcommand{\dx}[1]{\mathd #1}

\renewcommand{\acnst}[1]{}

\begin{document}

\title{An Orthogonality Property of the Legendre Polynomials} 
\author{L. Bos$^{1}$, A. Narayan$^{2}$, N. Levenberg$^{3}$ and
F. Piazzon$^{4}$}

\maketitle

\footnotetext[1]{Dept. of Computer Science, University of
Verona, Italy\\
e-mail: leonardpeter.bos@univr.it}

\footnotetext[2]{Mathematics Dept., University of Massachusetts Dartmouth, North Dartmouth, Massachusetts, USA\\
e-mail: akil.narayan@umassd.edu}

\footnotetext[3]{Dept. of Mathematics, Indiana University,
Bloomington, Indiana, USA\\
e-mail: nlevenbe@indiana.edu}

\footnotetext[4]{Dept. of Mathematics, University of
Padua, Italy\\
e-mail: fpiazzon@math.unipd.it}

\begin{abstract}
  We give a remarkable additional othogonality property of the classical Legendre polynomials on the real interval $[-1,1]$: polynomials up to degree $n$ from this family are mutually orthogonal under the arcsine measure weighted by the degree-$n$ normalized Christoffel function.
\end{abstract}

\vskip0.2cm
\noindent
{\small{\bf Keywords:}
Legendre polynomials, Christoffel function, equilbrium measure\/}

\vskip0.2cm
\noindent
{\small{\bf AMS Subject Classification:}
33C45, 41A10, 65C05

\vskip1cm
Let $P_n(x)$ denote the classical Legendre polynomial of degree $n$ and
\[P_n^*(x):=\frac{\sqrt{2n+1}}{\sqrt{2}}P_n(x)\]
its orthonormalized version. 
Thus, with $\delta_{i,j}$ the Kronecker delta, the family $P^\ast_n$ satisfies
  {
  \begin{align*}
    \int_{-1}^1 P^\ast_i(x) P^\ast_j(x)\, \dx{x} &= \delta_{i,j}, & i, j &\geq 0
  \end{align*}
}

We consider the normalized (reciprocal of) the 
 associated Christoffel function
 \begin{equation}\label{CF}
 K_n(x):=\frac{1}{n+1}\sum_{k=0}^n (P_k^*(x))^2.
 \end{equation}
 As is well known, $K_n(x)dx$ tends weak-$*$ to the equilibrium measure of complex potential theory for the interval $[-1,1],$ and more precisely,
 \[\lim_{n\to\infty}K_n(x)=\frac{1}{\pi}\frac{1}{\sqrt{1-x^2}}, \quad x\in(-1,1)\]
 locally uniformly. In other words,
 \[\lim_{n\to\infty}\frac{1}{K_n(x)}\frac{1}{\pi}\frac{1}{\sqrt{1-x^2}}=1, \quad x\in(-1,1)\]
locally uniformly, and it would not be  unexpected that, at least asymptotically,
\[\int_{-1}^1P_i^*(x)P_j^*(x) \frac{1}{K_n(x)}\frac{1}{\pi}\frac{1}{\sqrt{1-x^2}}dx
\approx \delta_{ij},\quad 0\le i,j\le n.\]
The purpose of this note is to show that the above is actually an identity.
\begin{theorem}\label{thm:main-result} With the above notation
\begin{equation}\label{MR}
\int_{-1}^1P_i^*(x)P_j^*(x) \frac{1}{K_n(x)}\frac{1}{\pi}\frac{1}{\sqrt{1-x^2}}dx= \delta_{ij},\quad 0\le i,j\le n.
\end{equation}
\end{theorem}
We expect this result to have use in applied approximation problems. For example, one application lies in polynomial approximation of functions from point-evaluations. Our result indicates that the functions $\left\{Q_j(x)\right\}_{j=0}^n \triangleq \left\{ \frac{1}{\sqrt{K_n(x)}} P^\ast_j(x) \right\}_{j=0}^n$ are an orthonormal set on $(-1,1)$ under the Lebesgue density $\frac{1}{\pi\sqrt{1-x^2}}$. If we generate Monte Carlo samples from this density, evaluate an unknown function at these samples, and perform least-squares regression using $Q_j$ as a basis, then a stability factor for this problem is given by $\max_{x \in [-1,1]} \sum_{j=0}^n Q_j^2(x) = n+1$ \cite{cohen_stability_2013}. In fact, this is the smallest attainable stability factor, and therefore this approximation strategy has optimal stability. 

The remainder of this document is devoted to the proof of \eqref{MR}.

\noindent {\bf Proof of Theorem \ref{thm:main-result}.} We change variables letting $x=\cos(\theta)$ to arrive at the equivalent expression
\[\frac{1}{\pi}\int_{0}^\pi\frac{P_i^*(\cos(\theta))P_j^*(\cos(\theta))} {K_n(\cos(\theta))}d\theta= \delta_{ij},\quad 0\le i,j\le n\]
which by symmetry is equivalent to
\begin{equation}\label{InTheta}
\frac{1}{2\pi}\int_{0}^{2\pi}\frac{P_i^*(\cos(\theta))P_j^*(\cos(\theta))} {K_n(\cos(\theta))}d\theta= \delta_{ij},\quad 0\le i,j\le n.
\end{equation}
Now, for $z\in\C,$ let $\displaystyle{J(z):=(z+1/z)/2}.$ Then for $z=e^{i\theta}$ in the integral \eqref{InTheta} we obtain, $\displaystyle{d\theta =-iz^{-1}dz},$
$\cos(\theta)=J(z),$ and the equation becomes
\begin{equation}\label{Inz}
\frac{1}{2\pi i}\int_C\frac{z^{-1}P_i^*(J(z))P_j^*(J(z))} {K_n(J(z))}d z= \delta_{ij},\quad 0\le i,j\le n
\end{equation}
where $C$ is the unit circle, oriented in the counter-clockwise direction.

The proof is a direct calculation of \eqref{Inz} based on the following lemmas.

First note that $K_n(\cos(\theta))$ is a {\it positive} trigonometric polynomial (of degree $2n$). By the F\'ejer-Riesz Factorization Theorem there exists a trigonometric polynomial, $T_n(\theta)$ say, such that
\[K_n(\cos(\theta))=|T_n(\theta)|^2.\]
In general the coefficients of the factor polynomial, $T_n(\theta)$ in this case, are algebraic functions of the coefficients of the original polynomial. However in our case we have the explict (essentially) rational factorization.

\begin{prop} \label{FR} (F\'ejer-Riesz Factorization of $K_n(J(z))$) Let
\begin{equation}\label{Fn}
F_n(z):=\frac{d}{dz}\left(z^{n+1}P_n(J(z))\right)=(n+1)z^nP_n(J(z))+\frac{z^{n-1}(z^2-1)}{2}P_n'(J(z)).
\end{equation}
Then 
\begin{equation}\label{KnFact}
K_n(J(z))=\frac{1}{2(n+1)}F_n(z)F_n(1/z).
\end{equation}
\end{prop}
\noindent {\bf Proof.} To begin, one may easily verify that
\begin{equation}\label{Fn(1/z)}
F_n(1/z)=z^{-2n}\{(n+1)z^nP_n(J(z))-\frac{z^{n-1}(z^2-1)}{2}P_n'(J(z))\}.
\end{equation}
Hence
\begin{align*}
F_n(z)F_n(1/z)&=z^{-2n}\{(n+1)^2z^{2n}(P_n(J(z)))^2-z^{2(n-1)}\left(
\frac{z^2-1}{2}\right)^2(P_n'(J(z))^2\}\\
&=(n+1)^2(P_n(J(z)))^2-z^{-2}\left(
\frac{z^2-1}{2}\right)^2(P_n'(J(z))^2.
\end{align*}
Now notice that
\begin{align*}
z^{-2}\left(\frac{z^2-1}{2}\right)^2 &= \frac{1}{4} \left(z - \frac{1}{z}\right)^2\\
&= \frac{1}{4} \left(z^2 + 2 + \frac{1}{z^2} - 4\right) \\
&= J^2(z) - 1
\end{align*}

so that
\[ F_n(z)F_n(1/z)=(n+1)^2(P_n(J(z)))^2-(J(z)^2-1)(P_n'(J(z)))^2.\]
The result follows then from  Lemma \ref{KnIdentity}, below. $\square$ 

\bigskip
\begin{lemma}\label{KnIdentity} For all $x\in\C,$ we have
\[K_n(x)=\frac{1}{2(n+1)}\left((n+1)^2(P_n(x))^2-(x^2-1)(P_n'(x))^2\right).\]
\end{lemma}
 \noindent{\bf Proof.} 
First, we collect the following known identities concerning Legendre polynomials \cite{szego_orthogonal_1939}: 
   \begin{subequations}
   \begin{align}
     \nonumber\allowdisplaybreaks[1]
     \shortintertext{{\small (Christoffel-Darboux formula)}}
     \label{eq:leg-cd} \sum_{k=0}^n P_k^2(x) &= \frac{n+1}{2} \left[ P_{n+1}'(x) P_n(x) - P_{n+1}(x) P_n'(x) \right] \\
     \shortintertext{{\small (Three-term recurrence)}}
     \label{eq:leg-rec} (n+1) P_{n+1}(x) &= (2 n + 1) x P_n(x) - n P_{n-1}(x) \\
     \shortintertext{{\small (Differentiated three-term recurrence)}}
     \label{eq:leg-drec} (n+1) P'_{n+1}(x) &= (2 n + 1) \left(P_n(x) + x P_n'(x)\right) - n P'_{n-1}(x) \\
     \label{eq:leg-christ} (x^2-1) P_n'(x)   &= n \left( x P_n(x)-P_{n-1}(x) \right) \\
     \label{eq:leg-ichrist} (2 n + 1) P_n(x) &= P_{n+1}'(x) - P_{n-1}'(x)
   \end{align}
 \end{subequations}

We easily see from the Christoffel-Darboux formula that
 \[K_n(x)=\frac{1}{2}\left(P_{n+1}'(x)P_n(x)-P_{n+1}(x)P_n'(x)\right).\]
Hence the result holds iff { \setstretch{0.5}
\begin{multline*}
(n+1)P_{n+1}'(x)P_n(x)-(n+1)P_{n+1}(x)P_n'(x) \\
= (n+1)^2(P_n(x))^2-(x^2-1)(P_n'(x))^2
\end{multline*}
\begin{center}$\Updownarrow$ \eqref{eq:leg-drec},\eqref{eq:leg-rec},\eqref{eq:leg-christ} \end{center} 
\begin{multline*}
  \left\{ (2 n + 1) \left(P_n(x) + x P_n'(x)\right) - n P'_{n-1}(x)\right\} P_n(x) 
  \\ - \left\{ (2 n + 1) x P_n(x) - n P_{n-1}(x) \right\} P_n'(x) \\
  =(n+1)^2(P_n(x))^2- n \left( x P_n(x)-P_{n-1}(x) \right) P_n'(x)
\end{multline*}
\begin{center}$\Updownarrow$ \end{center} 
\begin{multline*}
(2 n + 1) \left(P_n(x)\right)^2 - n P_{n-1}'(x) P_n(x) + n P_{n-1}(x) P_n'(x)  \\
  =(n+1)^2(P_n(x))^2- n x P_n(x) P_n'(x) + n P_{n-1}(x) P_n'(x)
\end{multline*}
\begin{center}$\Updownarrow$ \end{center} 
\begin{align*}
-n^2 \left(P_n(x)\right)^2 - n P_{n-1}'(x) P_n(x) = - n x P_n(x) P_n'(x) 
\end{align*}
\begin{center}$\Updownarrow$ \end{center} 
\begin{align*}
x P_n'(x)  = n P_n(x) + P_{n-1}'(x)
\end{align*}
\begin{center}$\Updownarrow$ \eqref{eq:leg-drec}\end{center} 
\begin{align*}
  \frac{1}{2 n + 1} \left( (n+1) P_{n+1}'(x) + n P_{n-1}'(x)\right) - P_n(x) = n P_n(x) + P_{n-1}'(x)
\end{align*}
\begin{center}$\Updownarrow$ \end{center} 
\begin{align*}
  (n+1) P_{n+1}'(x) = (n+1) P_{n-1}'(x) + (2n + 1) (n+1) P_n(x),
\end{align*}
\setstretch{1.0}
and this last relation is the same as the relation \eqref{eq:leg-ichrist}.
} $\square$

\bigskip
There is somewhat more that can be said about $F_n(z).$

\begin{lemma} Let $F_n(z)$ be the polynomial of degree $2n$ defined in \eqref{Fn}. Then all of its zeros are simple and lie in the {\bf interior} of the unit disk.
\end{lemma}
\noindent {\bf Proof.}
The polynomial 
\[Q_n(z):=z^{n+1}P_n(J(z))=z\left\{z^nP_n(J(z))\right\}\]
has a zero at $z=0$ and its other zeros are those of $P_n(J(z))$ which are those
$z\in \C$ for which $J(z)=r\in(-1,1),$ a zero of $P_n(x).$ But
\begin{align*}
&\quad J(z)=r\in (-1,1)\\
\iff&\quad (z+1/z)/2=r\\
\iff&\quad z^2-2rz+1=0\\
\iff&\quad z=r\pm i\sqrt{1-r^2}.
\end{align*} 
In particular $|z|=1$ for the zeros of $z^nP_n(J(z)).$ It follows then from
the Gauss-Lucas Theorem that the zeros of $F_n(z)$ are in the convex hull of $z=0$ and certain points on the unit circle, i.e., are all in the {\it closed} unit disk.

Suppose a zero of $F_n(z)$ satisfies $|z| = 1$. By Proposition \ref{FR}, 
\[K_n(J(z))=\frac{1}{2(n+1)}F_n(z)F_n(1/z),\]
so that $K_n(J(z))$ also vanishes. But $|z|=1$ implies that $J(z) \in [-1,1]$, and $K_n(J(z))$ thus cannot vanish. Therefore, no zeros of $F_n$ lie on the unit circle.

To see that the zeros are all simple, an elemenary calculation and the  ODE for Legendre polynomials gives us
\[F_n'(z)=2n(n+1)z^{n-1}P_n(J(z))+\{nz^n-(n+1)z^{n-2}\}P_n'(J(z)).\]
Hence $F_n(z)=F_n'(z)=0$ if and only if
\[\left(\begin{array}{cc}n+1&\frac{z^2-1}{2}\cr
\frac{2n(n+1)}{z}&nz-\frac{n+1}{z}\end{array}\right)
\left(\begin{array}{c}z^nP_n(J(z))\cr z^{n-1}P_n'(J(z))\end{array}\right)
=\left(\begin{array}{c}0\cr0\end{array}\right).\]
But the determinant of this matrix is
\[(n+1)(nz-(n+1)/z)-n(n+1)(z^2-1)/z=-(n+1)/z\neq0.\]
Hence $F_n(z)=F_n'(z)=0$ if and only if $z^nP_n(J(z))=z^{n-1}P_n'(J(z))=0$ if and only if $P_n(J(z))=P_n'(J(z))=0,$ which is not possible as $P_n(x)$ has only simple zeros. $\square$

\bigskip

The integral \eqref{Inz} can therefore be expressed as
\begin{align*}
\frac{1}{2\pi i}\int_C\frac{z^{-1}P_i^*(J(z))P_j^*(J(z))} {K_n(J(z))}d z&= 
\frac{1}{2\pi i}\int_C\frac{2(n+1)z^{-1}P_i^*(J(z))P_j^*(J(z))} {F_n(z)F_n(1/z)}dz\\
&=\frac{1}{2\pi i}\int_C\frac{2(n+1)z^{2n-1}P_i^*(J(z))P_j^*(J(z))} {F_n(z)z^{2n}F_n(1/z)}dz\\
&=\frac{1}{2\pi i}\int_C\frac{2(n+1)z^{2n-1}P_i^*(J(z))P_j^*(J(z))} {F_n(z)G_n(z)}dz
\end{align*}
where we define the {\it polynomial} of degree $2n,$
\begin{equation}\label{Gn}
G_n(z):=z^{2n}F_n(1/z).
\end{equation}
As all the zeros of $F_n(z)$ are in the interior of the unit disk, the zeros of 
$G_n(z)$ are all exterior to the (closed) unit disk.

The following formulas for $F_n(z)$ and $G_n(z)$ will be useful.

\begin{lemma}\label{FnGn}
We have
\[F_n(z)=\frac{z^n}{z^2-1}\{((2n+1)z^2-1)P_n(J(z))-2nzP_{n-1}(J(z))\}\]
and
\[G_n(z)=\frac{z^n}{z^2-1}\{(z^2-(2n+1))P_n(J(z))+2nzP_{n-1}(J(z))\}.\]
\end{lemma}
\noindent {\bf Proof.} From the formula \eqref{Fn} we have
\[F_n(z)=(n+1)z^nP_n(J(z))+z^{n-1}\frac{z^2-1}{2}P_n'(J(z))\]
and from \eqref{Fn(1/z)},
\[G_n(z)=(n+1)z^nP_n(J(z))-z^{n-1}\frac{z^2-1}{2}P_n'(J(z)).\]

{ From the Legendre Polynomial identity \eqref{eq:leg-christ} with $x = J(z)$, we obtain}
\[z^{n-1}\frac{z^2-1}{2}P_n'(J(z))=2n\frac{z^{n+1}}{z^2-1}
J(z)P_n(J(z))-
{2 n \frac{z^{n+1}}{z^2-1}}
P_{n-1}(J(z)).\]
Combining these gives the result.  $\square$

\bigskip

It is also interesting to note that $F_n(z)$ is a certain Hypergeometric function.

\begin{lemma} \label{PropsFn} We have 
\begin{enumerate}
\item The polynomial $y=F_n(z)$ is a solution of the ODE
\[(1-z^2)y''+2\frac{(n-2)z^2-n}{z}y'+6ny=0.\]
\item If $F_n(z)=:f_n(z^2)$ then $y=f_n(z)$ is a solution of the
Hypergeometric ODE
\[z(1-z)y''+(c-(a+b+1)z)y'-aby=0\]
with $a=-n,$ $b=3/2$ and $c=1/2-n.$
\item $\displaystyle{f_n(z)=2^{-2n}{2n \choose n}\, _2F_1(a,b;c;z).}$
\item $\displaystyle{F_n(z)=2^{-2n}{2n \choose n}\, _2F_1(a,b;c;z^2).}$
\item \begin{align*}
F_n(z)&=2^{-2n}{2n \choose n}\sum_{k=0}^n
\frac{(2k+1){n\choose k}^2}{{2n\choose 2k}} z^{2k}\\
&=2^{-2n}\sum_{k=0}^n(2k+1)
{2k \choose k}{2n-2k\choose n-k}z^{2k}.
\end{align*}
\item If we write $F_n(z)=\sum_{k=0}^nc_kz^{2k},$ then
\begin{align*}
G_n(z)&=\sum_{k=0}^n c_{n-k}z^{2k}\\
&=\sum_{k=0}^n \frac{2(n-k)+1}{2k+1}c_kz^{2k}.
\end{align*}
\end{enumerate}
\end{lemma}
\noindent {\bf Proof.} Equation (1) is easily verified using the ODE for $P_n(x)$ and the definition of $F_n(z),$ \eqref{Fn}. 

(2) follows by changing variables $z'=z^2.$

(3) follows as $_2F_1(a,b;c;z)$ is the only polynomial solution of the Hypergeometric ODE. The constant of proportionality is calculated by noting that the
leading coefficient of $_2F_1(a,b;c;z)$ is $2n+1$ whereas that of $f_n(z)$ is 
$(2n+1)/2^n$ times the leading coefficient of $P_n(x),$ i.e.,
$(2n+1)/2^n\times {2n \choose n}/2^n.$

(4) is trivial {from the defintion $f_n(z^2) := F_n(z)$.}
 
(5) follows from the fact that
\[_2F_1(a,b;c;z)=\sum_{k=0}^n \frac{(a)_k(b)_k}{(c)_k}\frac{z^k}{k!},\]
{with $(a)_k$ the rising factorial,} and calculating
\[ (a)_k=(-1)^k\frac{{n!}}{(n-k)!},\,\,(b)_k=2^{-2k}\frac{(2k+1)!}{k!},\,\,
(c)_k=(-1)^k2^{-2k}\frac{(2n)!(n-k)!}{(2n-2k)!n!}\]
so that
\[\frac{1}{k!}\frac{(a)_k(b)_k}{(c)_k}=\frac{(2k+1){n\choose k}^2}{{2n\choose 2k}}.\]

(6) follows easy from the fact that $G_n(z):=z^nF_n(1/z)=\sum_{k=0}^n
c_{n-k}z^k$ and that $c_{n-k}$ is easily computed from the formula for $c_k$ given in (5). $\square$

Returning the the proof of the Theorem, we will actually show that
\begin{description}
\item[(a)] $\displaystyle{\frac{1}{2\pi i}\int_C\frac{2(n+1)z^{2n-1}P_k(J(z))} {F_n(z)G_n(z)}dz=0,\quad 1\le k\le 2n,}$
\item[(b)] $\displaystyle{\frac{1}{2\pi i}\int_C\frac{2(n+1)z^{2n-1}P_0(J(z))} {F_n(z)G_n(z)}dz=2}.$
\end{description}
The Theorem follows directly as, for $i,j\le n$ we may expand
\[P_i^*(x)P_j^*(x)=\sum_{k=0}^{2n}a_kP_k(x)\]
for certain coefficients $a_k.$ From (a) and (b) we then conclude that
\[\frac{1}{2\pi i}\int_C\frac{2(n+1)z^{2n-1}P_i^*(J(z))P_j^*(J(z))} {F_n(z)G_n(z)}dz=2a_0.\]
But for $i\neq j,$
\[a_0=\frac{1}{2}\int_{-1}^1 P_i^*(x)P_j^*(x) P_0(x)dx=0\]
as $P_0(x)=1$ and $P_i^*(x)$ and $P_j^*(x)$ are orthogonal.
While for $i=j,$ we have
\[a_0=\frac{1}{2}\int_{-1}^1 (P_i^*(x))^2P_0(x)dx=\frac{1}{2}.\]

We will now compute the partial fraction decomposition of {the integrands in (a) and (b),}
\[\frac{2(n+1)z^{2n-1}P_k(J(z))} {F_n(z)G_n(z)},\]
which will involve the following two pairs of families of functions.

\begin{defn} \label{ABCD} We define
\[A_0^{(n)}(z):=z^{n-1},\quad A_1^{(n)}(z):=z^{n-2},\]
\[B_0^{(n)}(z):=z^{n-1},\quad B_1^{(n)}(z):=z^{n},\]
with
\[(n+k+1)A_{k+1}^{(n)}:=(2(n+k)+1)J(z)A_k^{(n)}(z)-(n+k)A_{k-1}^{(n)}(z),
\quad k=1,2,\cdots\]
and
\[(n+k+1)B_{k+1}^{(n)}:=(2(n+k)+1)J(z)B_k^{(n)}(z)-(n+k)B_{k-1}^{(n)}(z),
\quad k=1,2,\cdots.\]
Further, we let
\[C_0^{(n)}(z):=z^{n-1}, \quad C_1^{(n)}(z):=z^{n-2}\frac{(2n+1)z^2-1}{2n},\]
\[D_0^{(n)}(z):=z^{n-1}, \quad D_1^{(n)}(z):=z^{n-2}\frac{(2n+1)-z^2}{2n},\]
with
\[(n-k)C_{k+1}^{(n)}(z):=(2(n-k)+1)J(z)C_k^{(n)}(z)-(n-k+1)C_{k-1}^{(n)}(z), \quad 1\le k\le n-1\]
and
\[(n-k)D_{k+1}^{(n)}(z):=(2(n-k)+1)J(z)D_k^{(n)}(z)-(n-k+1)D_{k-1}^{(n)}(z),
\quad 1\le k\le n-1.\]
\end{defn}

\begin{prop}\label{PartFrac} We have
\begin{equation}\label{Pn+k}
2(n+1)z^{2n-1}P_{n+k}(J(z))=A_k^{(n)}(z)G_n(z)+B_k^{(n)}(z)F_n(z),\quad k=0,1,2,\cdots
\end{equation}
so that 
\[\frac{2(n+1)z^{2n-1}P_{n+k}(J(z))} {F_n(z)G_n(z)}=\frac{A_k^{(n)}(z)}{F_n(z)}
+\frac{B_k^{(n)}(z)}{G_n(z)}, \quad k=0,1,2,\cdots \]
and
\begin{equation}\label{Pn-k}
2(n+1)z^{2n-1}P_{n-k}(J(z))=C_k^{(n)}(z)G_n(z)+D_k^{(n)}(z)F_n(z),\quad 0\le k\le n
\end{equation}
so that
\[\frac{2(n+1)z^{2n-1}P_{n-k}(J(z))} {F_n(z)G_n(z)}=\frac{C_k^{(n)}(z)}{F_n(z)}
+\frac{D_k^{(n)}(z)}{G_n(z)},\quad 0\le k\le n.\]
\end{prop}
\noindent {\bf Proof} (by induction). As $A_0^{(n)}=C_0^{(n)}=z^{n-1}$ and $B_0^{(n)}=D_0^{(n)}=z^{n-1},$ the $k=0$ case is in common and we calculate
\begin{align*}
&=z^{n-1}(F_n(z)+G_n(z))\\
&{\stackrel{\textrm{{(Lemma \ref{FnGn})}}}{=}}z^{n-1}\frac{z^n}{z^2-1}\{((2n+1)z^2-1)+(z^2-(2n+1))\}P_n(J(z))\\
&=z^{2n-1}\frac{1}{z^2-1}\{2(n+1)z^2-2(n+1)\}P_n(J(z))\\
&=2(n+1)z^{2n-1}P_n(J(z)),
\end{align*}
as desired.

We now prove \eqref{Pn+k} for the case $k=1.$ We calculate,  using Lemma \ref{FnGn},
\begin{align*}
&A_1^{(n)}(z)G_n(z)+B_1^{(n)}(z)F_n(z)\\
&=z^{n-2}G_n(z)+z^nF_n(z)\\
&=z^{n-2}(G_n(z)+z^2F_n(z))\\
&=z^{n-2}\frac{z^n}{z^2-1}\{[(z^2-(2n+1))P_n(J(z))+2nzP_{n-1}(J(z))]\\
&\qquad + z^2[((2n+1)z^2-1)P_n(J(z))-2nzP_{n-1}(J(z))]\}\\
&=\frac{z^{2n-2}}{z^2-1}\{(2n+1)(z^4-1)P_n(J(z))-2nz(z^2-1)P_{n-1}(J(z))\}\\
&=z^{2n-2}\{(2n+1)(z^2+1)P_n(J(z))-2nzP_{n-1}(J(z))\}\\
&=z^{2n-1}\{(2n+1)(z+1/z)P_n(J(z))-2nP_{n-1}(J(z))\}\\
&=z^{2n-1}\{(2n+1)2J(z)P_n(J(z))-2nP_{n-1}(J(z))\}\\
&{\stackrel{\eqref{eq:leg-rec}}{=}}2(n+1)z^{2n-1}P_{n+1}(J(z)){.}
\end{align*}
Now for \eqref{Pn-k} for $k=1.$ We calculate, again using Lemma \ref{FnGn},
\begin{align*}
&C_1^{(n)}(z)G_n(z)+D_1^{(n)}(z)F_n(z)\\
&=z^{n-2}\left\{\frac{(2n+1)z^2-1}{2n}G_n(z)+\frac{(2n+1)-z^2}{2n}F_n(z)
\right\}\\
&=z^{n-2}\frac{z^n}{z^2-1}\left\{\frac{(2n+1)z^2-1}{2n}[(z^2-(2n+1))P_n(J(z))+2nzP_{n-1}(J(z))]\right.\\
&\qquad \left. + \frac{(2n+1)-z^2}{2n}[((2n+1)z^2-1)P_n(J(z))-2nzP_{n-1}(J(z))]\right\}\\
&=\frac{z^{2n-2}}{z^2-1}2nz\left\{\frac{(2n+1)z^2-1}{2n}-\frac{(2n+1)-z^2}{2n}\right\}P_{n-1}(J(z))\\
&={\frac{z^{2n-1}}{z^2-1}}((2(n+1)z^2-2(n+1))P_{n-1}(J(z))\\
&=2(n+1)z^{2n-1}P_{n-1}(J(z)).
\end{align*}

The rest of the proof proceeds by induction. Assuming that \eqref{Pn+k} and
\eqref{Pn-k} hold from 0 up to a  certain $k,$ we prove that they also hold for $k+1.$ To this
end we calculate
\begin{align*}
&A_{k+1}^{(n)}(z)G_n(z)+B_{k+1}^{(n)}(z)F_n(z)\\
&=\frac{(2(n+k)+1)J(z)A_k^{(n)}(z)-(n+k)A_{k-1}^{(n)}(z)}
{n+k+1}G_n(z)\\
&\quad +\frac{(2(n+k)+1)J(z)B_k^{(n)}(z)-(n+k)B_{k-1}^{(n)}(z)}{n+k+1}F_n(z)\\
&=\frac{1}{n+k+1}\{(2(n+k)+1)J(z)[A_k^{(n)}(z)G_n(z)+B_k^{(n)}(z)F_n(J(z))]\\
&\quad -(n+k)[A_{k-1}^{(n)}(z)G_n(z)+B_{k-1}^{(n)}F_n(J(z))]\\
&=2(n+1)z^{2n-1}\frac{1}{n+k+1}\{(2(n+k)+1)J(z)P_{n+k}(J(z))\\
&\quad -(n+k)P_{n+k-1}(J(z)){\}}\quad \hbox{(by the induction hypothesis)}\\
&=2(n+1)z^{2n-1}P_{n+k+1}(J(z)),
\end{align*}
by the three-term recursion formula for Legendre polynomials \eqref{eq:leg-rec} with degree $m=n+k.$

Similarly,
\begin{align*}
&C_{k+1}^{(n)}(z)G_n(z)+D_{k+1}^{(n)}(z)F_n(z)\\
&=\frac{(2(n-k)+1)J(z)C_k^{(n)}(z)-(n-k+1)C_{k-1}^{(n)}(z)}
{n-k}G_n(z)\\
&\quad +\frac{(2(n-k)+1)J(z)D_k^{(n)}(z)-(n-k+1)D_{k-1}^{(n)}(z)}{n-k}F_n(z)\\
&=\frac{1}{n-k}\{(2(n-k)+1)J(z)[C_k^{(n)}(z)G_n(z)+D_k^{(n)}(z)F_n(J(z))]\\
&\quad -(n-k+1)[C_{k-1}^{(n)}(z)G_n(z)+D_{k-1}^{(n)}F_n(J(z))]\\
&=2(n+1)z^{2n-1}\frac{1}{n-k}\{(2(n-k)+1)J(z)P_{n-k}(J(z))\\
&\quad -(n-k+1)P_{n-(k-1)}(J(z)){\}}\quad \hbox{(by the induction hypothesis)}\\
&=2(n+1)z^{2n-1}P_{n-(k+1)}(J(z)),
\end{align*}
using the reverse three-term recursion,
\[mP_{m-1}(x)=(2m+1)xP_m(x)-(m+1)P_{m+1}(x)\]
with $m=n-k.$ $\square$

\bigskip Due to the $J(z)$ factor in the recursive definitions of Definition \ref{ABCD}, the functions $A_k^{(n)}(z),$ $B_k^{(n)}(z),$ $C_k^{(n)}(z)$
and $D_k^{(n)}(z)$  are all Laurent polynomials. It is easy to verify that they have
the forms:
\begin{itemize}
  \item $\displaystyle{A_k^{(n)}(z)=\sum_{j={n-(k+1)}}^{{n+(k-3)}}a_kz^k,}$ $k\ge 1$
  \item $\displaystyle{B_k^{(n)}(z)=\sum_{j={n-(k-1)}}^{{n+(k-1)}}b_kz^k,}$ $k\ge 1$
  \item $\displaystyle{C_k^{(n)}(z)=\sum_{j={n-(k+1)}}^{{n+(k-1)}}c_kz^k,}$ $1\le k\le n$
  \item $\displaystyle{D_k^{(n)}(z)=\sum_{j={n-(k+1)}}^{{n+(k-1)}}d_kz^k,}$ $1\le k\le n.$
\end{itemize}
In particular, for $0\le k\le n-1$ the functions $A_k^{(n)}(z),$ $B_k^{(n)}(z),$ $C_k^{(n)}(z)$ are all {\it polynomials} of degree at most $2n-2.$

Hence, for $1\le k\le 2n-1,$
\[
\frac{1}{2\pi i}\int_C \frac{2(n+1)z^{2n-1}P_k(J(z))}{F_n(z)G_n(z)}dz=\frac{1}{2\pi i}
\left\{\int_C \frac{p(z)}{F_n(z)}dz + \int_C \frac{q(z)}{G_n(z)}dz\right\}\]
for certain polynomials $p(z)$ and $q(z)$ of degree at most $2n-2.$

Now, $\displaystyle{\int_C \frac{q(z)}{G_n(z)}dz=0}$ as all the zeros of
$G_n(z)$ lie outside the (closed) unit disk. 
Further, if we let $z_j$ $1\le j\le 2n,$ be the (simple) zeros of $F_n(z),$ we may write
\[\frac{p(z)}{F_n(z)}=\frac{p(z)/c_n}{F_n(z)/c_n}=\sum_{j=1}^{2n} \frac{R_j}{z-z_j}\]
where $c_n$ is the leading coefficient of $F_n(z).$ Hence
\[\frac{1}{2\pi i}\int C \frac{A_n(z)}{F_n(z)}dz=\frac{1}{2\pi i}2\pi i
\sum_{j=1}^{2n} R_j.\]
But $\sum_{j=1}^{2n} R_j$ is the leading coefficient (of $z^{2n-1}$) of $p(z)/c_n,$ i.e., $0,$ as $p(z)$ is of degree at most $2n-2$.  It follows that
\[\frac{1}{2\pi i}\int_C \frac{2(n+1)z^{2n-1}P_k(J(z))}{F_n(z)G_n(z)}dz=0,\quad 1\le k \le 2n-1.\]

The cases $P_0(J(z))$ and $P_{2n}(J(z))$ are special.

\medskip  First consider the case $P_{2n}(J(z)).$ We have from Proposition \ref{PartFrac}, with $k=n,$ 
\[2(n+1)z^{2n-1}P_{2n}(J(z))=A_n^{(n)}(z)G_n(z)+B_n^{(n)}(z)F_n(z).\]
However, $\displaystyle{A_n^{(n)}(z)=\sum_{j=-1}^{2n-3}a_jz^j}$ has a $z^{-1}$ while $\displaystyle{B_n^{(n)}(z)=\sum_{j=+1}^{2n-1}b_jz^j}$ is still a polynomial,
of degree at most $2n-1.$ Therefore it is still the case that
$\displaystyle{\int_C \frac{B_n^{(n)}(z)}{G_n(z)}dz=0}$  (the zeros of $G_n(z)$ being all oustide the unit disk). We need to show that $\displaystyle{\int_C \frac{A_n^{(n)}(z)}{F_n(z)}dz=0}.$ \acnst{But this is not difficult.} Write
$A_n^{(n)}(z)=q(z)+c/z$ where $q(z)$ is a polynomial of degree $2n-3.$
Then
\[\int_C \frac{A_n^{(n)}(z)}{F_n(z)}dz=\int_C \frac{q(z)}{F_n(z)}dz +c \int_C \frac{1}{zF_n(z)}dz.\]
The first integral on the right is zero as the coefficient in $q(z)$ of $z^{2n-1}$ is $0.$ For the second integral, decompose
\begin{align*}
c\int_C \frac{1}{zF_n(z)}dz&=c\int_C \frac{1}{F_n(0)}\left\{\frac{1}{z}-
\frac{(F_n(z)-F_n(0))/z}{F_n(z)}\right\}dz\\
&=\frac{1}{F_n(0)}\left\{\int_C \frac{1}{z}dz-\int_C \frac{(F_n(z)-F_n(0))/z}{F_n(z)}\right\}dz.
\end{align*}
The first integral is trivially $2\pi i,$ while the second is
$2\pi i$ times the coefficient of $z^{2n-1}$ in $(F_n(z)-F_n(0))/z$ divided by the
leading coefficient (of $z^{2n}$) in $F_n(z),$ i.e., $2\pi i \times 1.$ Hence, indeed
\[\int_C \frac{A_n^{(n)}(z)}{F_n(z)}dz=0.\]

\medskip

Lastly we calculate
\[\frac{1}{2\pi i}\int_C \frac{2(n+1)z^{2n-1}P_{0}(J(z))}{F_n(z)G_n(z)}dz=\frac{1}{2\pi i}\int_C \frac{2(n+1)z^{2n-1}}{F_n(z)G_n(z)}dz.\]
We still have  {\eqref{Pn-k} from Proposition \ref{PartFrac}},
\[2(n+1)z^{2n-1}P_{n-n}(z)=C_n^{(n)}(z)G_n(z)+D_n^{(n)}(z)F_n(z)\]
However,  $C_n^{(n)}(z)$ and $D_n^{(n)}(z)$ have the form
\[C_n^{(n)}(z)=\sum_{j=-1}^{2n-1}c_kz^k,\quad D_n^{(n)}(z)=\sum_{j=-1}^{2n-1}d_kz^k,\]
i.e., are both of the form $p(z)+c/z$ where $p(z)$ is a polynomial of degree
$2n-1.$  Specifically, write
$C_n^{(n)}(z)=p(z)+c/z$ and $D_n^{(n)}(z)=q(z)+d/z.$ {We thus have the following expression:
  \begin{align*}
    \frac{1}{2\pi i} \int_C \frac{2 (n+1)z^{2n-1} P_0(z)}{F_n(z) G_n(z)} \dx{z} &= 
    \frac{1}{2\pi i} \int_C \frac{C_n^{(n)}(z)}{F_n(z)} \dx{z} + \frac{1}{2 \pi i} \int_C \frac{D_n^{(n)}(z)}{G_n(z)} \dx{z} \\
    &= \frac{1}{2\pi i} \left[ \int_C \frac{p(z)}{F_n(z)} \dx{z} + \int_C \frac{c}{z F_n(z)} \dx{z} \right. \\
    & \hskip 25pt \left. + \int_C \frac{q(z)}{G_n(z)} \dx{z} + \int_C \frac{d}{z G_n(z)} \dx{z} \right].
  \end{align*}
}
{We need to show that this expression has value 2.} Now we have already shown that $\displaystyle{\int_C \frac{1}{zF_n(z)}dz=0}$ and also remarked that
$\displaystyle{\int_C \frac{q(z)}{G_n(z)}dz=0}.$ Hence we need to calculate
$\displaystyle{\frac{1}{2\pi i}\int_C \frac{p(z)}{F_n(z)}dz}$
and $\displaystyle{\frac{1}{2\pi i}\int_C \frac{d}{zG_n(z)}dz}.$
But the first of these is just the (leading) coefficient of $z^{2n-1}$ in $p(z),$ i.e.,
in $C_n^{(n)}(z)$ divided by the leading coefficient of $F_n(z).$ From Lemma 
\ref{PropsFn} we have that
\[F_n(z)=2^{-2n}(2n+1){2n \choose n}z^{2n}+\cdots\]
and
it is easy to verify by induction that also
$\displaystyle{C_n^{(n)}(z)=2^{-2n}(2n+1){2n \choose n}z^{2n}+\cdots}.$
Hence
\[\frac{1}{2\pi i}\int_C \frac{p(z)}{F_n(z)}dz=1.\]
For the second integral, decompose as before
\[
\frac{1}{2\pi i}\int_C \frac{d}{zG_n(z)}dz=\frac{d}{2\pi i}\frac{1}{G_n(0)}
\int_C\left\{\frac{1}{z}-\frac{(G_n(z)-G_n(0))/z}{G_n(z)}\right\}dz.
\]
The second integral above is zero as all the zeros of $G_n(z)$ are outside the closed unit disk. The first integral is trivially $d/G_n(0).$ From Lemma \ref{PropsFn} parts 4. and 6. we have that 
\[G_n(0)=(2n+1)F_n(0)=2^{-2n}(2n+1){2n \choose n}\]
whereas the coefficient of $z^{-1}$ in $B_n^{(n)}(z)$ is easily verified by induction
to have the same value. Hence
\[\frac{1}{2\pi i}\int_C \frac{d}{zG_n(z)}dz=1\]
and we have shown that
\[\frac{1}{2\pi i}\int_C \frac{2(n+1)z^{2n-1}P_{0}(J(z))}{F_n(z)G_n(z)}dz=1+1=2,\]
as claimed. $\square$

\bibliographystyle{amsplain}
\bibliography{LegendreOrthogonality}

\providecommand{\bysame}{\leavevmode\hbox to3em{\hrulefill}\thinspace}
\providecommand{\MR}{\relax\ifhmode\unskip\space\fi MR }
\providecommand{\MRhref}[2]{%
  \href{http://www.ams.org/mathscinet-getitem?mr=#1}{#2}
}
\providecommand{\href}[2]{#2}
\begin{thebibliography}{1}

\bibitem{cohen_stability_2013}
Albert Cohen, Mark~A. Davenport, and Dany Leviatan, \emph{On the {Stability}
  and {Accuracy} of {Least} {Squares} {Approximations}}, Foundations of
  Computational Mathematics \textbf{13} (2013), no.~5, 819--834.

\bibitem{szego_orthogonal_1939}
G\'{a}bor Szeg\"{o}, \emph{Orthogonal {Polynomials}}, 4th ed., American
  Mathematical Soc., 1939.

\end{thebibliography}

\end{document}